\documentclass[11pt]{article}

\input{amssym.def}
\input{amssym}
\usepackage{eufrak}
\usepackage{color}

\newtheorem{theorem}{Theorem}

\newtheorem{proposition}[theorem]{Proposition}
\newtheorem{remark}[theorem]{Remark}

\def\pa{{\partial}}

\def\Y{\mathbf{Y}}
\def\YR{\mathbf{Y}(R)}

\def\qq{q^{-1}}

\def\hLL{{\widehat{\cal L}}}
\def\hL{\hat L}
\def\LL{{\cal{L}}}
\def\Lc{L_{can}}

\def\A{{\cal A}}

\def\ot{\otimes}

\def\C{{\Bbb C}}

\def\End{{\rm End}}
\def\vv{V^{\otimes 2}}

\def\hQ{\hat Q}

\def\al{{\alpha}}

\def\be{\begin{equation}}
\def\ee{\end{equation}}

\begin{document}

\makeatletter
\renewcommand{\theequation}{{\thesection}.{\arabic{equation}}}
\@addtoreset{equation}{section} \makeatother

\title{Drinfeld-Sokolov  reduction  in  quantum algebras}

\author{\rule{0pt}{7mm} Dimitri
Gurevich\thanks{gurevich@ihes.fr}\\
{\small\it LAMAV, Universit\'e de Valenciennes,
59313 Valenciennes, France}\\
\rule{0pt}{7mm} Pavel Saponov\thanks{Pavel.Saponov@ihep.ru}\\
{\small\it
National Research University Higher School of Economics,}\\
{\small\it 20 Myasnitskaya Ulitsa, Moscow 101000, Russian Federation}\\
{\small \it and}\\
{\small \it
Institute for High Energy Physics, NRC "Kurchatov Institute"}\\
{\small \it Protvino 142281, Russian Federation}\\
\rule{0pt}{7mm} Dmitry Talalaev\thanks{dtalalaev@yandex.ru}\\
{\small\it
Moscow State University, Faculty of Mechanics and Mathematics}\\
{\small\it 119991 Moscow, Russian Federation}}

\maketitle

\begin{abstract}
Applying the method of the paper \cite{CT}, we perform a quantum version of the Drinfeld-Sokolov reduction in Reflection Equation algebras and braided Yangians,
associated with involutive and Hecke symmetries of general forms. This reduction is  based on the  Cayley-Hamilton identity valid for the
generating matrices of these algebras.
\end{abstract}



\section{Introduction}
In the seminal paper \cite{DS} it was shown that any connection operator $\pa_u-M(u)$, where  $\pa_u=\frac{d}{du}$ and $M(u)$ is a $N\times N$ matrix smoothly depending of
the parameter $u$, can be reduced by means of the gauge transformations
$$
\pa_u-M(u)\,\mapsto\, g(u) (\pa_u-M(u)) g(u)^{-1}=\pa_u-g(u) M(u) g(u)^{-1}+g(u) \pa_u(g(u)^{-1})
$$
to the form $\pa_u-M_{can}(u)$, where $M_{can}(u)$  looks like (\ref{can}). This form of operators and matrices will be called {\em canonical} or, more precisely, {\em second canonical}.
Also, the authors of \cite{DS} found that the reduced Poisson structure can be identified with  the second Gelfand-Dickey one. In \cite{FF} this Poisson reduction was quantized in
terms of the BRST cohomology and the quantum object was identified with $W$-algebras.

A natural question arises: what are suitable  operators and their canonical forms related to the quantum groups or other quantum algebras. Thus, in \cite{FR, FRS} a version of the Drinfeld-Sokolov
reduction is considered, where the operator $\pa_u$ is replaced by that $D_q f(x)=f(qx)$, $q\not=1$ and the  gauge transformations are replaced by
$$
D_q-M(u)\,\mapsto \,g(qu)(D_q-M(u))g(u)^{-1}.
$$
However, this structure is not a deformation of the classical one since for $q=1$ the corresponding operators differ from the classical opertors.

In this note we apply another approach, suggested in \cite{CT}. There, the second canonical form was found for a matrix $M(u)$ subject to the relation
\be
\big[M_1(u), M_2(v)]=\left[\frac{P}{u-v}, M_1(u)+M_2(v)\right],
\label{Gaud}
\ee
where $P$ is the $N^2\times N^2$ matrix of the usual flip, $M(u)=\sum_{k\geq 0} M[k] u^{-k}$ is an $N\times N$ matrix, and, as usual, $M_1=M\ot I$, $M_2=I\ot M$. Note that the Hamiltonians
of the rational Gaudin model are expressed via such a matrix.

Besides, the matrix $M(u)$ satisfies the Cayley-Hamilton (CH) identity with modified matrix powers of $M(u)$. Namely, the terms $M(u)^k$ have been replaced by some combinations of this matrix
and its derivatives in $u$. By using these modified powers one can show that the operators $\pa_u-M(u)$ and $\pa_u-M_{can}(u)$ are in a sense similar (see \cite{CT} for detail). Namely, in
this sense the term "Drinfeld-Sokolov reduction" was used in \cite{CT}.

In \cite{CF} a version of the CH identity for the matrices $L(u)$ which generates the Drinfeld's Yangian $\Y(gl(N))$ was found (below, $L(u)$ and other similar matrices will be called {\em generating
matrices}). The powers of this matrix are also understood in a modified sense: namely, in the products of the matrix $L(u)$ the arguments of the factors are shifted. The use of the CH identity allows
one to define the second canonical form of the generating matrix $L(u)$ and to establish the similarity of the initial and reduced forms of this matrix.

The aim of this note is to define an analogous Drinfeld-Sokolov (DS) reduction for the generating matrices of some quantum algebras. Namely, we perform this reduction for the generating matrices
of the Reflection Equation algebras, corresponding to constant (i.e. independent of the parameters) involutive or Hecke symmetries $R$ and for the {\it braided Yangians}, also arising from these
symmetries. Such type Yangian-like algebras have been recently introduced \cite{GS2, GS3}.

We mention here one of possible applications of this result of particular interest. The quantum algebras under consideration have direct relevance to the theory of quantum integrable systems due to the presence of typical commutative families. We hope that the construction of the second canonical form will allow us to establish its relations with new versions of the KZ equation.

\medskip

\noindent
{\bf Acknowledgements:} The work of D.T. was partially supported by the grant of the Simons foundation and the RFBR grant 17-01-00366 A.
The work of P.S. has been funded by the Russian Academic Excellence Project '5-100' and was also partially supported by the RFBR grant 16-01-00562.

\section{Quantum Matrix algebras}
Let us briefly describe the quantum algebras under consideration. First, recall that by a current $R$-matrix one usually means an operator $R(u,v)$ depending on parameters $u$ and $v$ and
subject to the so-called Quantum Yang-Baxter equation
$$
R_{12}(u,v) R_{23}(u,w) R_{12}(v,w) =R_{23}(v,w) R_{12}(u,w) R_{23}(u,v),
$$
where $R_{12}(u,v)=R(u,v)\ot I$ and $R_{23}(u,v)=I\ot R(u,v)$. If $R$ is independent of the parameters it is also called {\em a braiding}. In this case we shall assume the operator $R:\vv\to\vv$
to be either involutive $R^2=I$ or to satisfy the Hecke condition
$$
(q\, I-R)( \qq\, I+R)=0,  \quad q\in \C, \quad q^2\not=1.
$$
Here, $V$ is a vector space of the dimension $N$ (over the field $\C$). These braidings are respectively called {\em involutive and  Hecke symmetries}.

The best known are the Hecke symmetries coming from the quantum groups $U_q(sl(N))$. However, there are known numerous examples of involutive and Hecke symmetries which are deformations
neither of usual nor of super-flips (see \cite{GPS2} and the references therein).

All symmetries, we are dealing with, are assumed to be {\em skew-invertible} and to have the bi-rank $(m|0)$.  For the definitions of the notions "skew-invertible" and "bi-rank" the reader is referred to
\cite{GS2, GS3}. We want only to mention that for a skew-invertible braiding $R$, an $R$-trace
$$
{\rm Tr}_R:\End(V)\to \C
$$
in a sense coordinated with $R$ can be defined. Moreover, for any matrix $A$ the expression ${\rm Tr}_R \, A$ is also well-defined. The properties of such an $R$-traces can be found in \cite{O, GPS2}.

Let us respectively define two quantum matrix algebras by the following systems of  relations on their generators
\be
R\, T_1\,T_2-T_1\, T_2\, R=0,\quad T=\|t_i^j\|_{1\leq i,j \leq N}
\label{RTT}
\ee
\be
R\, L_1\, R\,L_1-L_1\, R\, L_1\, R=0,\quad L=\|l_i^j\|_{1\leq i,j \leq N}
\label{RE}
\ee
where $R$ is assumed to be a skew-invertible involutive or Hecke symmetry. The former algebra is called an RTT one, the latter one is called the Reflection Equation (RE) algebra.
Their detailed consideration can be found in \cite{GPS2}. Here, we want only to observe that if $R$ is a deformation of the usual flip $P$, the dimensions of the homogenous  components of both
algebras are classical, i.e. equal to those in ${\rm Sym}(gl(N))$ (if $R$ is a Hecke symmetry, the value of $q$ is assumed to be generic).

Now, let us exhibit Yangian-like algebras associated with current $R$-matrices. First, observe that the current 
$R$-matrices, we are dealing with, are constructed from involutive or Hecke symmetries via the so-called Baxterization procedure described in \cite{GS2, GS3}. 
This procedure results in the following $R$-matrix
\be
R(u,v) = R - g(u,v) I,\quad {\rm where}\quad g(u,v)=\frac{1}{u-v}\quad {\rm or}\quad g(u,v)=\frac{(q-q^{-1})u}{u-v}.
\label{curR}
\ee
If $R$ is an involutive symmetry, $g(u,v)$ is defined by the former formula. If $R$ is a Hecke symmetry, $g(u,v)$ is defined by the latter one.

The corresponding Yangian-like algebras, introduced in \cite{GS2, GS3} are also of two types. They are respectively defined by the following systems
\be
R(u,v)\, T_1(u)\, T_2(v)=T_1(v)\, T_2(u)\, R(u,v),
\label{ya}
\ee
\be
R(u,v)\, L_1(u)\, R\, L_1(v)=L_1(v)\, R\, L_1(u)\, R(u,v),
\label{brya}
 \ee
where the generating matrices $T(u)$ and $L(u)$ are assumed to possess a series expansion as the matrix $M(u)$ above. We call the former (resp., latter) algebras {\em the Yangians of
RTT type} (resp., {\em the braided Yangians}). Note that the braided Yangians are defined similarly to the RE algebras, but with the current $R$-matrices at the outside positions.

Let us introduce the following notation
\be
L_{\overline 1}(u)=L_1(u),\quad L_{\overline k}(u)=R_{k-1} L_{\overline{k-1}}(u) R^{-1}_{k-1},\quad k\geq 2,
\label{note}
\ee
where we write $R_i$ instead of $R_{i\, i+1}$ (recall that  $R_{i\, i+1}$ stands for the operator $R$ acting at the $i$-th and $i+1$-th positions in the tensor product $V^{\ot k}$, $i+1\leq k$).
In the RE algebras we use the same notation for the generating matrix $L$ which is independent of the parameters.

By using this notation we can present the defining relations of an RE algebra in the form
$$
R_1 \, L_{\overline 1} \, L_{\overline 2}- L_{\overline 1} \, L_{\overline 2}\, R_1=0.
$$
Moreover, in this algebra the following holds:
$$
R_k \, L_{\overline k} \, L_{\overline{k+1}}- L_{\overline k} \, L_{\overline{k+1}}\, R_k=0, \quad \forall\,k\ge 1.
$$

This notation enables us to define analogs of symmetric polynomials in the RTT and RE algebras in a uniform way.
Thus, the {\em power sums} are respectively defined as follows
$$
p_k(T)={\rm Tr}_{R(12\dots k)} R_{k-1}R_{k-2}\dots R_2R_1T_1\,T_2\dots T_k,
$$
\be
p_k(L)={\rm Tr}_{R(12\dots k)} R_{k-1}R_{k-2}\dots  R_2R_1 L_{\overline{1}}L_{\overline 2}\dots L_{\overline{k}}.
\label{ps}
\ee
Here the notation ${\rm Tr}_{R(12\dots k)}$ means that the $R$-traces are applied at the positions $1,2, \dots ,k$. Note that in the RE algebra the formula (\ref{ps}) can be simplified
to $p_k(L)={\rm Tr}_R L^k$ whereas for the power sums $p_k(T)$ in the RTT algebra such a transformation is not possible. 

In a similar manner the ``quantum powers" of the matrices $T$ and $L$ can be defined:
$$
T^{[k]}:={\rm Tr}_{R(2\dots k)} R_{k-1}R_{k-2}R_2\dots R_1T_1T_2\dots T_k,
$$
$$
L^{[k]}:={\rm Tr}_{R(2\dots k)} R_{k-1}R_{k-2}\dots  R_2R_1 L_{\overline{1}}L_{\overline{2}}\dots L_{\overline{k}}.
$$
However, if the former formula cannot be simplified, the latter one can be reduced to the usual one: $L^{[k]}=L^k$.

Also, exhibit analogs of elementary symmetric polynomials in both algebras
$$
e_0(T)=1,\qquad  e_k(T):={\rm Tr}_{R(1\dots k)} (\A^{(k)} T_{1} \,T_{ 2}\dots T_{ k}),\quad k\geq 1,
$$
\be
e_0(L)=1,\qquad  e_k(L):={\rm Tr}_{R(1\dots k)} (\A^{(k)} L_{\overline 1} \,L_{\overline 2}\dots L_{\overline k}),\quad k\geq 1.
\label{ele}
\ee
Here  $\A^{(k)}: V^{\ot k} \to V^{\ot k}$, $k\ge 1$ are the skew-symmetrizers (i.e. the projectors of skew-symmetriza\-tion)
which are recursively defined by
\be
\A^{(1)} = I,\quad \A^{(k)} =\frac{1}{k_q} \,\A^{(k-1)}\left(q^{k-1} I-(k-1)_q R_{k-1}\right)\A^{(k-1)},\quad k\geq 2,
\label{q-asym}
\ee
Hereafter, we use the standard notation: $k_q=\frac{q^k-q^{-k}}{q-q^{-1}}$.

Note that if the bi-rank of $R$ is $(m|0)$, $m\geq 2$, the skew-symmetrizers $\A^{(k)}$ are trivial for $k> m$ and the rank of the skew-symmetrizer $\A^{(m)}$ is equal to 1.

Analogs of matrix powers, power sums and elementary symmetric polynomials can be also defined for generalized Yangians of both classes (see \cite{GS2, GS3}). We exhibit
them below for the braided Yangians. The point is that only for the braided Yangians and RE algebras we can use the CH identity for performing the DS reduction
in a way similar to that of \cite{CT}. This is due to the fact that in these algebras the powers of the generating matrices are classical (up to shifts of the arguments in the generating
matrices) and the structure of the CH identity is close to the classical one.

Observe that  analogs of the CH identity for the generating matrices of quantum matrix algebras associated with couples of braidings were presented in \cite{IOP}. Note also, that
quantum analogs of the elementary symmetric polynomials and power sums can be defined in the frameworks of the so-called half-quantum algebras as defined in \cite{IO}. However,
for proving their commutativity it is necessary to impose more strong relations giving rise to one of the quantum matrix algebras.

\section{DS Reduction in RE algebras}

Let $R$ be again a skew-invertible involutive or Hecke symmetry. Denote $\LL(R)$ the RE (\ref{RE}). As was shown in \cite{GPS1}, the generating matrix $L$ of the algebra $\LL(R)$
meets the {\em quantum CH identity} $Q(L)=0$,  where the {\it characteristic polynomial} $Q(t)$ reads
\be
Q(t)=t^m-q t^{m-1} e_1(L)+q^2 t^{m-2}e_2(L)+\dots +(-q)^{m-1} t e_{m-1}(L)+(-q)^{m}  e_m(L)=0.
\label{char}
\ee
Here, the factors $e_k(L)$ are the quantum elementary symmetric polynomials defined by (\ref{ele}).

\begin{remark}\rm
We stress a very important property of the polynomial  (\ref{char}): its coefficients belong to the center $Z(\LL(R))$ of the algebra $\LL(R)$. Let us introduce``eigenvalues" $\{\mu_i\}_{1\le i\le m}$
of the matrix $L$  in a natural way
$$
e_1(L) = \mu_1+...+\mu_m,\quad\dots\quad e_m(L) = \mu_1\cdot\dots \cdot \mu_m.
$$
These "eigenvalues"  are  elements of an algebraic extension of the center   $Z(\LL(R))$.  Consider  the quotient algebra
$$
\LL(R)/\langle e_1(L)-\al_1,\dots,e_m(L)-\al_m \rangle,\qquad \al_i\in \C,
$$
where $<I>$ stands for  the ideal generated by a set $I\subset \LL(R)$. This quotient is a  quantum analog of an orbit (or a union of orbits) in the coadjoint representation of the group $GL(N)$.
In \cite{GS1} there was considered the problem: for which values of $\al_i$ this quotient is an analog of a regular orbit. If it is so, we introduce the diagonal matrix $diag(\mu_1,\dots,\mu_m)$, where
the elements $\mu_i$ solve the system
 $e_1(L)=\al_1,\dots,e_m(L)=\al_m$,
 and treat this matrix as the {\em first canonical form}  of the generating matrix $L$.
\end{remark}

Let us define the {\em second canonical form} of the matrix $L$:
\be
L_{can}=\left(\begin{array}{ccccc}
0&1&0&...&0\\
0&0&1&...&0\\
...&...&...&...&...\\
0&0&0&...&1\\
a_m& a_{m-1}  &a_{m-2}&...&a_1
\end{array}\right) , \label{can}
\ee
where
$$
a_k=-(-q)^k e_k(L).
$$

Following \cite{CT} we show that the matrices $L$ and $\Lc$ are in a sense similar. Let $v\in V$ be an arbitrary non-trivial vector written as a one-row matrix. Then we introduce
the following $N\times N$ matrix
\be C=\left(\begin{array}{c}
v\\
vL\\
...\\
vL^{m-1}\end{array}\right). \label{CC}
\ee

\begin{proposition}
The following relation holds true
\be
C\,L=\Lc \, C. \label{simi}
\ee
\end{proposition}

In this sense we say that the matrices $L$ and $\Lc$ are {\em similar}.

\begin{remark}\rm
In order to justify the term "similar" it would be desirable to show that at least for some vectors $v$ the matrix $C$ is invertible in the skew-field of the algebra $\LL(R)$.
A similar problem is also open for other quantum algebras considered below.  \end{remark}

Along with the RE algebra $\LL(R)$ define its quadratic-linear deformation $\hLL(R)$ by the following system
\be
R\, \hL_1\, R\,\hL_1-\hL_1\, R\, \hL_1\, R=R\,\hL_1-\hL_1\, R,\qquad \hL=\|\hat{l}_i^j\|_{1\leq i,j \leq N}.
\label{mRE}
\ee

If the algebra $\LL(R)$ is a braided analog of the algebra ${\rm Sym}(gl(N))$, then $\hLL(R)$ is a braided analog of the enveloping algebra $U(gl(N))$. Namely, if $R$ is a
Hecke symmetry which is a deformation of the usual flip (for instance, that coming  from the quantum group $U_q(sl(N))$), the algebras $\LL(R)$ and $\hLL(R)$ turn into ${\rm Sym}(gl(N))$
and  $U(gl(N))$ respectively as $q\to 1$. By using the fact that the algebras
$\LL(R)$ and $\hLL(R)$ are isomorphic to each other if $q\not= \pm 1$, it is possible to get a characteristic polynomial  $\hQ(t)$ for the generating matrix
$\hL$ of the algebra $\hLL(R)$. Namely, we have (see \cite{GS2, GS3})
$$
{\hat Q}(t)= {\rm Tr}_{R(1\dots m)} \left(\A^{(m)}(tI-\hL_{\overline 1})((q^2t-q) I-\hL_{\overline 2})\dots ((q^{2(m-1)} t-q^{m-1}(m-1)_q) I-\hL_{\overline m})\right).
$$

\begin{proposition}
In the algebra $\hLL$ the following matrix identity takes place ${\hat Q}(L)=0$.
\end{proposition}

Passing to the limit $q\to 1$ we get the characteristic polynomial for the generating matrix\footnote{Note, that the defining relations on the generators $m_i^j$ of the algebra
$U(gl(N))$
$$
m_i^j\, m_k^l- m_k^l\,  m_i^j =m_i^l\delta_k^j- m_k^j\delta_i^l,
$$
can be written in a matrix form with the use of the generating $N\times N$ matrix $M=\|m_i^j\|$:
$$
 P\, M_1\, P\,M_1-M_1\, P\, M_1\, P=P\,M_1-M_1\, P.
$$}
 $M$ of the algebra $U(gl(N))$ (here $N=m$)
$$
{\cal Q}(t)= {\rm Tr}_{(1\dots N)} \left(\A^{(N)}(tI-M_1)((t-1) I-M_2)\dots ((t-N+1) I-M_N)\right),
$$
where $\A^{(N)}$ is the usual skew-symmetrizer in $V^{\otimes N}$ and ${\rm Tr}$ is the usual trace.

\begin{proposition}
In the algebra $U(gl(N))$ the following holds ${\cal Q}(M)=0$.
\end{proposition}

Note that  the famous Capelli determinant is defined by a similar formula. The same claim is valid for any algebra $\hLL$ provided  $R$ be an involutive symmetry, which can be
approximated by a Hecke symmetry.

Besides, it is possible to introduce the second canonical forms $\hL_{can}$ and $M_{can}$ of the matrices $\hL$ and $M$ respectively generating the algebras $\hLL(R)$ and $U(gl(N))$,
and to perform a DS reduction of the matrices $\hL$ and $M$. Upon replacing the matrix $L$ in  (\ref{CC}) by  $\hL$ and $M$ respectively, we get formulae similar to (\ref{simi}). However,
in the case of the algebra $\hLL$ we have first to normalize the polynomial ${\hat Q}(t)$ in order to get a monic one. Thus, if this monic polynomial is $L^m+b_{1} L^{m-1}+\dots + b_m I$,
we have to put $a_k=-b_k$ in the matrix (\ref{can}).

\begin{remark} \rm
There are other matrices with entries from the algebras under consideration for which the CH identities exist. First, consider the enveloping algebra $U(gl(N))$. Its generating matrix $M$
belongs to the so-called Kirillov's quantum family algebra (see \cite{K})
\be
(U(gl(N))\ot \End(V))^{gl(N)}.
\label{fam}
\ee
Upon replacing $V$ by other irreducible $U(gl(N))$-modules, we get other quantum family algebras. For their generating matrices a characteristic polynomial can be found. Note
that  in the algebra $\hLL(R)$ a q-analog of (\ref{fam}) is
$$
(\hLL\ot \End(V))^{U_q(sl(N))}.
$$
Changing $\hLL$ for $\LL$, we get  a ``$q$-family algebra" for the algebra $\LL(R)$.
 \end{remark}

\section{DS reduction in braided Yangians}

Now, consider a braided Yangians  $\YR$. Using the notations (\ref{note}) we can represent the defining system of this algebra as follows
$$
R(u,v)\, L_{\overline 1}(u)\,  L_{\overline  2}(v)=L_{\overline 1}(v)\,  L_{\overline 2}(u)\, R(u,v).
$$

First, we assume $R$ to be an involutive symmetry and the corresponding $R$-matrix $R(u,v)$ to be given by the first formula (\ref{curR}).

In this case we get  the recurrence formula defining the skew-symmetrizers ${\cal A}^{(k)}$  by putting $q=1$ in  (\ref{q-asym}). Thus, we have
$$
{\cal A}^{(k)} =\frac{1}{k} \, {\cal A}^{(k-1)}\left( I-(k-1)  R_{k-1}\right) {\cal A}^{(k-1)}.
$$

Let us respectively introduce  the {\em quantum elementary symmetric} elements and the {\em quantum matrix powers} of $L(u)$ as follows:
$$e_0(u)=1,\quad   e_k(u) = {\rm Tr}_{  R(1\dots k)}\left(  {\cal A}^{(k)} L_{\overline 1}(u)L_{\overline 2}(u-1)\dots L_{\overline k}(u-(k-1))\right),\quad k\geq 1,
$$
$$
L^{[0]}(u)=I,\quad   L^{[k]}(u) = L(u-(k-1))\dots L(u-1)L(u),\quad  k\geq 1.
$$
Then the following CH identity is valid for the generating matrix $L(u)$ (see \cite{GS2, GS3})
\be
\sum_{k=0}^m (-1)^k L^{[m-k]}(u-k)   e_k(u)=0. \label{CH1}
\ee

Consider the operator $L(u)e^{\pa_u}$, where $\pa_u=\frac{d}{du}$. Our current aim is to give an explicit canonical form $L_{can}(u)e^{\pa_u}$ of this operator and to show
that the operators $L(u)e^{\pa_u}$ and $L_{can}(u)e^{\pa_u}$ are similar in the sense of formula (\ref{simi}).

Note, that as $L_{can}(u)$ we use the matrix transposed to (\ref{can}). This is motivated by the fact that the coefficients $e_p(u)$ in the CH identity (\ref{CH1}) are on the right
hand side from the factors $L^{[k]}(u)$. (Note that in the RE algebras these canonical forms are equivalent due to the centrality of the coefficients $\sigma_k(L)$.)

So, the matrix $L_{can}(u)$ reads
\be
L_{can}(u)=\left(\begin{array}{cccccc}
0&0&...&0&0&a_m(u)\\
1&0&...&0&0&a_{m-1}(u)\\
0&1&...&0&0&a_{m-2}(u)\\
...&...&...&...&...&...\\
0&0&...&1&0&a_2(u)\\
0&0&...&0&1&a_1(u)
\end{array}\right) . \label{can1}
\ee

Let again $v\in V$ be an arbitrary non-trivial vector. Constitute an $N\times N$ matrix
$$
C(u)=\left(v, L(u)v, L^{[2]}(u+1)v\dots L^{[k]}(u+k-1))v\dots       L^{[m-1]}(u+m-2)v \right)
$$
where  $v$ stands for the corresponding column.

\begin{proposition} The following holds
\be
L(u)e^{\pa_u} C(u)= C(u) L_{can}(u)e^{\pa_u},
\label{simi1}
\ee
provided the entries $a_k(u)$ are
$$
a_k(u)=(-1)^{k+1} \, e_k(u+m-1).
$$
\end{proposition}

\noindent
{\bf Proof.} By using  the  evident relation
$$
e^{\pa_u}f(u)=f( u+1)e^{\pa_u},
$$
and canceling the operator $e^{\pa_u}$, we can present the relation (\ref{simi1}) as
\be
L(u)C( u+1)= C(u) L_{can}(u).
\label{simi2}
\ee

The equality of the corresponding columns from (\ref{simi2}), except the last ones, immediately follows from the relation
\be
L(u)L^{[k]}(u+k)=L^{[k+1]}(u+k).
\label{add}
\ee
This of the last columns follows from the CH identity  (\ref{CH1}).\hfill\rule{6.5pt}{6.5pt}
\medskip

Now, assume $R(u,v)$ to be defined by the second formula (\ref{curR}). Let us respectively define the {\em quantum elementary symmetric} elements and {\em quantum matrix powers} by
$$
e_0(u)=1,\quad e_k(u) = {\rm Tr}_{R(1\dots k)}\left(\A^{(k)} L_{\overline 1}(u)L_{\overline 2}(q^{-2}u)\dots L_{\overline k}(q^{-2(k-1)}u)\right),\quad k\geq 1,
$$
$$
L^{[0]}(u)=I,\quad L^{[k]}(u)= L(q^{-2(k-1)} u) L(q^{-2(k-2)} u)\dots L(u), \quad k\geq 1.
$$

Then the following form of the  CH identity is valid  for the generating matrix $L(u)$ (see \cite{GS2, GS3})
\be
\sum_{k=0}^m (-q)^k L^{[m-k]}(q^{-2k}u) e_k(u)=0.
\label{CH2}
\ee

Consider the operator $L(u)q^{2u  \pa_u}$ and constitute the following matrix
 \be
C(u)=\left(v, L(u)v, L^{[2]}(q^2u)v\dots L^{[k]}(q^{2(k-1)}u)v \dots  L^{[m-1]}(q^{2(m-2)}u)v \right),
\ee
where $v$ stands for the corresponding column.

\begin{proposition}
The following identity holds
\be
L(u)q^{2u\pa_u} C(u)= C(u) L_{can}(u)q^{2u\pa_u},
\label{simi3}
\ee
where the matrix $L_{can}(u)$ is given by (\ref{can1}) provided the entries $a_k(u)$ are defined by
$$
a_k(u)=-(-q)^{k} \, e_k(q^{2(m-1)}u).
$$
\end{proposition}

\noindent
{\bf Proof.} Now, we use the following  relation
$$
q^{2u\pa_u}f(u)=f(q^2 u)q^{2u\pa_u}.
$$
Then, upon canceling the factor $q^{2u\pa_u}$, we rewrite the relation (\ref{simi3}) as
\be
L(u)C(q^{2}u)= C(u) L_{can}(u).
\label{simi4}
\ee

Also, instead of (\ref{add}) we use the relation
$$
L(u)L^{[k]}(q^{2k} u)=L^{[k+1]}(q^{2k} u).
$$
This relation  entails  the equality of the corresponding columns from (\ref{simi4}), except the last ones. The equality of the last columns follows from
(\ref{CH2}). \hfill\rule{6.5pt}{6.5pt}

In conclusion, we want to remark that the above version of the DS reduction is valid in a  more general setting than the conventional  one.
 Indeed, it is not possible to reduce the generating matrix $M$ of the
enveloping algebra $U(gl(N))$ to a simpler  form by usual transformations $M\to g\, M\, g^{-1}$, whereas a version of the reduction based on the CH identity is valid.

\end{document}